\documentstyle[12pt]{article}

\makeatletter
\newdimen\normalarrayskip
\newdimen\minarrayskip
\normalarrayskip\baselineskip
\minarrayskip\jot
\newif\ifold \oldtrue 
\def\arraymode{\ifold\relax\else\displaystyle\fi}
\def\eqnumphantom{\phantom{(\theequation)}}
\def\@arrayskip{\ifold\baselineskip\z@\lineskip\z@
\else
\baselineskip\minarrayskip\lineskip2\minarrayskip\fi}
\def\@arrayclassz{\ifcase \@lastchclass \@acolampacol \or
\@ampacol \or \or \or \@addamp \or
\@acolampacol \or \@firstampfalse \@acol \fi
\edef\@preamble{\@preamble
\ifcase \@chnum
\hfil$\relax\arraymode\@sharp$\hfil
\or $\relax\arraymode\@sharp$\hfil
\or \hfil$\relax\arraymode\@sharp$\fi}}
\def\@array[#1]#2{\setbox\@arstrutbox=\hbox{\vrule
height\arraystretch \ht\strutbox
depth\arraystretch \dp\strutbox
width\z@}\@mkpream{#2}\edef\@preamble{\halign
\noexpand\@halignto
\bgroup \tabskip\z@ \@arstrut \@preamble \tabskip\z@ \cr}%
\let\@startpbox\@@startpbox \let\@endpbox\@@endpbox
\if #1t\vtop \else \if#1b\vbox \else \vcenter \fi\fi
\bgroup \let\par\relax
\let\@sharp##\let\protect\relax
\@arrayskip\@preamble}
\def\eqnarray{\stepcounter{equation}%
\let\@currentlabel=\theequation
\global\@eqnswtrue
\global\@eqcnt\z@
\tabskip\@centering
\let\\=\@eqncr
$$%
\halign to \displaywidth\bgroup
\eqnumphantom\@eqnsel\hskip\@centering
$\displaystyle \tabskip\z@ {##}$%
\global\@eqcnt\@ne \hskip 2\arraycolsep
$\displaystyle\arraymode{##}$\hfil
\global\@eqcnt\tw@ \hskip 2\arraycolsep
$\displaystyle\tabskip\z@{##}$\hfil
\tabskip\@centering
&{##}\tabskip\z@\cr}
\begingroup\ifx\undefined\newsymbol \else\def\input#1 {\endgroup}\fi

\def\nl#1#2{\mathop{#2}\limits_{#1}}

\newcommand{\pl}{\partial}

\newcommand{\ol}{\overline}
\newcommand{\bqq}{\begin{equation} \label}
\newcommand{\eeq}{\end{equation}}

\begin{document}

\setcounter{footnote}{1}
\def\thefootnote{\fnsymbol{footnote}}

\large
\centerline{\Large{\bf{Quadratic first integrals of geodesic equations}}}\smallskip
\centerline{\Large{\bf{of the special type 6-dimensional $h$-spaces }}}

\bigskip

\centerline{{\large
Zolfira Zakirova}\footnote{Kazan State University,
e-mail: Zolfira.Zakirova@soros.ksu.ru, Zolfira.Zakirova@ksu.ru}}

\bigskip

\abstract{\small  In this note we find a 6-dimensional
$h$-spaces of the $[(21\ldots1)(21\ldots1)\ldots(1\ldots1)]$ type and then determine quadratic first integrals
of the geodesic equations of these $h$-spaces.}

\begin{center}
\rule{5cm}{1pt}
\end{center}

\bigskip
\setcounter{footnote}{0}
\section{Introduction}

The aim of this paper is to investigate the $6$-dimensional pseudo-Riemannian space $V^6(g_{ij})$
with signature
$[++----]$. In particular, we find a metrics of 6-dimensional
$h$-spaces of the $[(21\ldots1)(21\ldots1)\ldots(1\ldots1)]$ type and then determine quadratic first integrals
of the geodesic equations of these $h$-spaces. The metrics of $h$-spaces of the $[2211]$
type has been obtained by the author \cite{zak1}.
The general method of determining pseudo-Riemannian manifolds that
admit some nonhomothetic projective group $G_r$ has been developed
by A.V.Aminova \cite{am2} and has got a detailed account in the paper
by the author \cite{zak} for the $h$-spaces of the type $[51]$.

In order to find the $h$-spaces admitting a nonhomothetic
infinitesimal projective transformation, one needs to integrate the
Eisenhart equation
\bqq{1}
h_{ij,k}=2g_{ij} \varphi_{,k}+g_{ik} \varphi_{,j}+
g_{jk} \varphi_{,i},
\eeq
which, in the skew-normal frames, is of the form \cite{am2}
\bqq{2}
X_r \overline a_{pq}+\sum_{h=1}^n e_h( \overline a_{hq} \gamma _{\tilde hpr}+
\overline a_{ph} \gamma_{\tilde  hqr})=\overline g_{pr}X_q \varphi+\overline
g_{qr}X_p \varphi
\quad
(p, q, r=1,\ldots,n),
\eeq
where
$$
X_r \varphi \equiv {\nl{r}\xi}^i \frac{\partial \varphi}
{\partial x^i},
\quad
\gamma_{pqr}=-\gamma_{qpr}=
{\nl{p}{\xi}\hspace{-2.5mm}\phantom{a}_{i,j}}{\nl{q}\xi}^i{\nl{r}\xi}^j,
\quad
a_{ij}=h_{ij}-2\varphi g_{ij},
$$
${\nl{i}\xi}^j$ are the components of skew-normal frames,
$\overline g_{pr}$ and
$\overline a_{pq}$ are the canonical forms of the tensors
$g_{pr}$, $a_{pq}$, respectively. For the type  $[(21\ldots1)(21\ldots1)\ldots(1\ldots1)]$
$h$-spaces these latter have the following form \cite{pet}
\bqq{3}
\ol{g}_{ij} dx^i dx^j=2 e_2 dx^1 dx^2+2 e_4 dx^3 dx^4 +
\sum_{\sigma=5}^6 e_{\sigma} {dx^{\sigma}}^2,
\eeq
$$
\ol{a}_{ij} dx^i dx^j=2 e_2 \lambda_2 dx^1 dx^2+e_2 {dx^2}^2+
2 e_4 \lambda_4 dx^3 dx^4+e_4 {dx^4}^2+
\sum_{\sigma=5}^6 e_{\sigma}  \lambda_{\sigma} {dx^{\sigma}}^2,
$$
$$
(e_1=e_2, e_3=e_4),
\quad
e_i =\pm 1,
$$
where $\lambda_1=\lambda_2$, $\lambda_3=\lambda_4$, $\lambda_5$, $\lambda_6$
are real functions, which may coincide. These functions are roots of the
characteristic equation ${\det}(h_{ij}-{\lambda}g_{ij})=0$.

\section{The metrics of the $h$-spaces of the $[22(11)]$ type}

Substituting the canonical forms $\overline g_{pr}$ and
$\overline a_{pq}$ from (\ref{3}) into (\ref{2}) and taking into account that
for the type $[22(11)]$  $h$-spaces $\lambda_5=\lambda_6$, $\tilde1=2$,  $\tilde2=1$,  $\tilde3=4$,
$\tilde4=3$,  $\tilde5=5$,  $\tilde6=6$, one gets the following
system of equations
\bqq{4}
X_r \lambda_2=0  \; (r \ne 2),
\quad
X_r \lambda_4=0  \; (r \ne 4),
\quad
X_r \lambda_6=0,
\eeq
$$
X_2(\lambda_2-\varphi)=X_4(\lambda_4-\varphi)=0,
\quad
\gamma_{121}=e_2 X_2 \phi,
\quad
\gamma_{343}=e_4 X_4 \phi,
$$
$$
\gamma_{142}=\gamma_{241}= \frac{e_2 X_4 {\varphi}}{\lambda_2-\lambda_4},
\quad
\gamma_{242}=-\frac{e_2 X_4 \varphi}{(\lambda_2-\lambda_4)^2},
\quad
\gamma_{324}=\gamma_{423}= \frac{e_4 X_2 {\varphi}}{\lambda_4-\lambda_2},
$$
$$
\gamma_{424}=-\frac{e_4 X_2 \varphi}{(\lambda_4-\lambda_2)^2},
\quad
\gamma_{244}=\frac{e_4 X_2 \varphi}{(\lambda_2-\lambda_4)^2},
\quad
\gamma_{s \sigma \sigma}=\frac{e_{\sigma}  X_s \varphi}{(\lambda_s-\lambda_{\sigma})},
$$
where $\sigma=5, 6$, $s=2, 4$, $\gamma_{56r}$ are arbitrary. All others $\gamma_{pqr}$ are equal to zero.

The system of linear partial differential equations
$$
X_q \theta={\nl{q}\xi}^i\pl_i \theta=0,
\quad
(q=1,\ldots,m, i=1,\ldots,6),
$$
where ${\nl{q}\xi}^i$ are the components of skew-normal frames, is
completely integrable if and only if the following conditions are fulfilled
(see \cite{ezen1} and \cite{am2})
\bqq{5}
(X_q, X_r) \theta=X_q X_r \theta-X_r X_q \theta=
\sum_{p=1}^6 e_p(\gamma_{pqr}-\gamma_{prq}) X_{\tilde p}.
\eeq
Using formulas (\ref{4}) and (\ref{5}), one can write down
the commutators of the
operators $X_i$ $(i=1,\ldots,6)$ in the $h$-spaces under consideration,
$$
(X_1,X_2)=-e_1\gamma_{121}X_2,
\quad
(X_1,X_3)=0,
\quad
(X_2,X_3)=e_4\gamma_{423}X_3,
$$
\bqq{6}
(X_1,X_4)=-e_2\gamma_{241}X_1,
\quad
(X_3,X_4)=-e_3\gamma_{343}X_4,
\eeq
$$
(X_2,X_4)=-e_2\gamma_{242}X_1-e_1\gamma_{142}X_2+e_4\gamma_{424}X_3+
e_3\gamma_{324}X_4,
$$
$$
(X_p,X_\sigma)=-e_{\tau} \gamma_{\tau \sigma p}X_{\tau},
\quad
(X_q,X_\sigma)=e_{\sigma} \gamma_{\sigma q \sigma}X_{\sigma}-
e_{\tau} \gamma_{\tau \sigma q}X_{\tau},
$$
$$
(X_5,X_6)=-e_5\gamma_{565}X_5+e_6\gamma_{656}X_6,
$$
where $p=1, 3$, $q=2, 4$, $\sigma, \tau=5, 6$ $(\sigma \ne \tau)$.

After coordinate transformation $x^{i'}=\theta^i(x)$, where $\theta^i$ are
solutions of the completely integrable systems from (\ref{6}),
one gets
\bqq{7}
{\nl{p}\xi}^i=P_p(x){\delta_p}^i,
\quad
{\nl{2}\xi}^{3}={\nl{2}\xi}^4={\nl{2}\xi}^{\sigma}={\nl{4}\xi}^1=
{\nl{4}\xi}^2={\nl{4}\xi}^{\sigma}={\nl{\sigma}\xi}^{\alpha}=0,
\eeq
where $\alpha=1, 2, 3, 4$.

Now integrating the system of equations (\ref{6}) with using
(\ref{4}) and (\ref{7}), one can find the components
${\nl{i}\xi}^j$ of skew-normal frames.
Using these results and the formulas \cite{am2}
\bqq{8}
g^{ij}=\sum_{h=1}^n e_h {\nl{h}{\xi}}^i {\nl{\tilde h}{\xi}}^j,
\quad
{\nl{h}{\xi}\hspace{-2.5mm}\phantom{a}_{i}}=g_{ij} {\nl{h}\xi}^j,
\quad
a_{ij}=\sum_{h, l=1}^n {e_h e_l} {\ol a}_{hl}
{\nl{\tilde h}{\xi}\hspace{-2.5mm}\phantom{a}_{i}} {\nl{\tilde l}{\xi}\hspace{-2.5mm}\phantom{a}_{j}}
\eeq
we can calculate the components of tensors $g_{ij}$ and $a_{ij}$.
Then we come to the following

\bigskip

\noindent
{\bf Theorem 1.} {\it If a symmetric tensor $h_{ij}$ of
the characteristics ${\rm [22(11)]}$ and a function $\varphi$
satisfy in $V^6(g_{ij})$
equation {\rm (\ref{1})}, then there exists a holonomic coordinate system
so that the function $\varphi$ and the tensors $g_{ij}$,
$h_{ij}$ are defined by the formulas
\bqq{10}
g_{ij}dx^idx^j=e_2A(f_4-f_2)\lbrace (f_4-f_2) dx^1 dx^2-
A (dx^2)^2)\rbrace+
\eeq
$$
e_4\tilde A(f_2-f_4)\lbrace (f_2-f_4) dx^3 dx^4-
\tilde A (dx^4)^2\rbrace+F_{\sigma \tau} (f_2-\lambda)^2(f_4-\lambda)^2
dx^{\sigma}dx^{\tau},
$$
\bqq{11}
a_{ij} dx^i dx^j=f_2 g_{i_1j_1}dx^{i_1}dx^{j_1}+
Ag_{12} (dx^2)^2+f_4 g_{i_2j_2}dx^{i_2}dx^{j_2}+
\tilde Ag_{34} (dx^4)^2+\lambda g_{\sigma \tau} dx^{\sigma}dx^{\tau},
\eeq
\bqq{12}
h_{ij}=a_{ij}+2\varphi g_{ij},
\quad
2\varphi=2f_2+2f_4+c,
\eeq
\bqq{13}
A=\epsilon x^1+\theta(x^2),
\quad
\tilde A=\tilde \epsilon x^3+\omega(x^4),
\eeq
where $\epsilon, \tilde \epsilon =0, 1$, $f_2=\epsilon x^2$, $f_4=\tilde \epsilon x^4+a$,
$\lambda$, $c$ and $a$ are some constants,
$a \ne 0$ when $\tilde{\epsilon}=0$,
$F_{\sigma \tau}(x^5, x^6), \theta(x^2), \omega(x^4)$ are arbitrary functions, $\theta \ne 0$ when
$\epsilon =0$, $\omega \ne 0$ when $\tilde \epsilon=0$, $i_1, j_1=1, 2$,
$i_2, j_2=3, 4$, $\sigma, \tau=5, 6$, $e_2, e_4=\pm 1$.}

\section{The metrics of the $h$-spaces of the $[2(21)1]$,\\
$[2(211)]$ types}

Substituting the canonical forms $\overline g_{pr}$ and
$\overline a_{pq}$ from (\ref{3}) into (\ref{2}) and taking into account that
for the type $[2(21)1]$  $h$-spaces $\lambda_4=\lambda_5$, one gets the following
system of equations
\bqq{13.1}
X_r \lambda_2=0  \; (r \ne 2),
\quad
X_r \lambda_5=0,
\quad
X_r \lambda_6=0 \; (r \ne 6),
\eeq
$$
X_2(\lambda_2-\varphi)=X_6(\lambda_6-\varphi)=0,
\quad
\gamma_{121}=e_2 X_2 \phi,
$$
$$
\gamma_{162}=\gamma_{261}= \frac{e_2 X_6 {\varphi}}{\lambda_2-\lambda_6},
\quad
\gamma_{262}=-\frac{e_2 X_6 \varphi}{(\lambda_2-\lambda_6)^2},
\quad
\gamma_{3s4}=\gamma_{4s3}= \frac{e_4 X_s {\varphi}}{\lambda_5-\lambda_s},
$$
$$
\gamma_{4s4}=-\frac{e_4 X_s \varphi}{(\lambda_5-\lambda_s)^2},
\quad
\gamma_{2\sigma \sigma}=\frac{e_{\sigma}  X_2 \varphi}{(\lambda_2-\lambda_{\sigma})},
\quad
\gamma_{565}=\frac{e_5 X_6 \varphi}{\lambda_5-\lambda_6},
$$
where  $\sigma=5, 6$, $s=2, 6$, $\gamma_{45r}$ are arbitrary,
all others $\gamma_{pqr}$ are equal to zero.

The commutators of the
operators $X_i$ $(i=1,\ldots,6)$ in the $h$-spaces of the  $[2(21)1]$ type have the form
\bqq{14}
(X_1,X_2)=-e_1\gamma_{121}X_2,
\quad
(X_1,X_3)=0,
\quad
(X_2,X_3)=e_4\gamma_{423}X_3,
\eeq
$$
(X_1,X_4)=-e_5\gamma_{541}X_5,
\quad
(X_1,X_5)=-e_4\gamma_{451}X_3,
\quad
(X_1,X_6)=-e_2\gamma_{261}X_1,
$$
$$
(X_2,X_4)=e_3\gamma_{324}X_4+e_4\gamma_{424}X_3-e_5\gamma_{542}X_5,
\quad
(X_2,X_5)=e_5\gamma_{525}X_5-e_4\gamma_{425}X_3,
$$
$$
(X_2,X_6)=-e_2\gamma_{262}X_1-e_1\gamma_{162}X_2+e_6\gamma_{626}X_6,
\quad
(X_3,X_4)=-e_5\gamma_{543}X_5,
$$
$$
(X_3,X_5)=-e_4\gamma_{453}X_3,
\quad
(X_3,X_6)=-e_4\gamma_{463}X_3,
\quad
(X_4,X_5)=-e_4\gamma_{454}X_3+e_5\gamma_{545}X_5,
$$
$$
(X_4,X_6)=-e_3\gamma_{364}X_4-e_4\gamma_{464}X_3+e_5\gamma_{546}X_5,
\quad
(X_5,X_6)=e_4\gamma_{456}X_3-e_5\gamma_{565}X_5.
$$
The systems $X_i \theta=0\;(i\ne2),
\;X_j \theta=0\;(j\ne 4),  \;X_k \theta=0\;(k\ne6)$ are completely integrable and have the following
solutions: $\theta^2$, $\theta^4$,  $\theta^6$.
The systems $X_3 \theta=X_4 \theta=X_5\theta=X_6\theta=0$,
$X_1 \theta=X_2 \theta=X_3\theta=X_6\theta=0$ and
$X_1\theta=X_2\theta=X_6\theta=0$  are completely integrable too. The first system has the solutions
$\theta^1$ and $\theta^2$,
the second system has the solutions $\theta^4$ and  $\theta^5$, the third system
has the solutions $\theta^3$, $\theta^4$ and $\theta^5$.
After coordinate transformation $x^{i'}=\theta^i(x)$, one gets
\bqq{15}
{\nl{p}\xi}^i=P_p(x){\delta_p}^i,
\quad
{\nl{2}\xi}^{s}={\nl{4}\xi}^q={\nl{5}\xi}^q={\nl{5}\xi}^4=0,
\eeq
where $p=1, 2, 3$, $s=3, 4, 5, 6$, $q=1, 2, 6$.

Now integrating the system of equations (\ref{14}) with using
(\ref{13.1}) and (\ref{15}), and then calculating the components of tensors $g_{ij}$ and $a_{ij}$
we come to the following results
\bqq{16}
g_{ij}dx^idx^j=e_2\lbrace 2(f_6-f_2)A dx^1 dx^2-
A^2 (dx^2)^2\rbrace+
\eeq
$$
(f_6-\lambda)(f_2-\lambda)^2\lbrace2e_4dx^3 dx^4-
e_4(\Sigma+\omega)(dx^4)^2+e_5 (dx^5)^2\rbrace+
e_6 (f_2-f_6)^2 (dx^6)^2,
$$
\bqq{17}
a_{ij} dx^i dx^j=f_2 g_{i_1j_1}dx^{i_1}dx^{j_1}+
g_{12} (dx^2)^2+\lambda g_{i_2j_2}dx^{i_2}dx^{j_2}+
g_{34} (dx^4)^2+f_6 g_{66} (dx^6)^2,
\eeq
\bqq{18}
h_{ij}=a_{ij}+(2f_2+f_6+c)g_{ij},
\quad
\varphi=f_2+\frac{1}{2}f_6+c,
\eeq
\bqq{19}
A=\epsilon x^1+\theta(x^2),
\quad
\Sigma=2(f_2-\lambda)^{-1}+(f_6-\lambda)^{-1},
\eeq
where $\epsilon=0, 1$, $f_2=\epsilon x^2$, $\lambda$ and $c$ are some
constants, $\theta(x^2), \omega(x^4, x^5), f_6(x^6)$ are arbitrary function, $\theta \ne 0$
when $\epsilon =0$, $i_1, j_1=1, 2$,
$i_2, j_2=3, 4, 5$, $e_2, e_4, e_5, e_6=\pm1$.

After same computations for $h$--spaces of the  $[2(211)]$ type, we obtain
\bqq{20}
g_{ij}dx^idx^j=2e_2A dx^1 dx^2+
(f_2-\lambda)^2\lbrace2e_4dx^3 dx^4-
e_4(\Sigma+\omega)(dx^4)^2+g_{\sigma \tau} dx^{\sigma} dx^{\tau}\rbrace,
\eeq
\bqq{21}
a_{ij} dx^i dx^j=2f_2 g_{12}dx^1dx^2+
g_{12} (dx^2)^2+\lambda g_{pq}dx^pdx^q+g_{34} (dx^4)^2,
\eeq
\bqq{22}
h_{ij}=a_{ij}+(2f_2+c)g_{ij},
\quad
\varphi=f_2+c,
\eeq
\bqq{23}
A=\epsilon x^1+\theta(x^2),
\quad
\Sigma=2(f_2-\lambda)^{-1},
\eeq
where $\epsilon=0, 1$, $f_2=\epsilon x^2$,
$\lambda, c$ are some constants,
$\theta(x^2), \omega(x^4, x^5, x^6), g_{\sigma \tau}(x^4, x^5, x^6)$ are arbitrary function, $\theta \ne 0$
when $\epsilon =0$, $p, q=3, 4, 5, 6$, $\sigma, \tau=5, 6$, $e_2, e_4=\pm1$.

Let us formulate these results in following theorem.

\bigskip

\noindent
{\bf Theorem 2.} {\it If a symmetric tensor $h_{ij}$ of
the characteristics $[2(21)1]$, $[2(211)]$ and a function $\varphi$
satisfy in $V^6(g_{ij})$ the Eisenhart equation, then there exists a holonomic coordinate system
so that the function $\varphi$ and the tensors $g_{ij}$,
$h_{ij}$ are defined by formulas {\rm (\ref{16})-(\ref{23})}.}

\section{The metrics of the $h$-spaces of the $[(22)11]$,\\
$[(221)1]$  types}

For the $h$--spaces of the $[(22)11]$ type from (\ref{2}), we get the following system of equations
$$
X_r \lambda_4=0,
\quad
X_r \lambda_{\sigma}=0 \; (r \ne \sigma),
\quad
X_{\sigma}(\lambda_{\sigma}-\varphi)=0,
$$
\bqq{24}
\gamma_{14r}=\gamma_{23r},
\quad
\gamma_{1\sigma 2}=\gamma_{2 \sigma 1}= \frac{e_2 X_{\sigma} {\varphi}}{\lambda_4-\lambda_{\sigma}},
\quad
\gamma_{3 \sigma 4}=\gamma_{4 \sigma 3}= \frac{e_4 X_{\sigma} {\varphi}}{\lambda_4-\lambda_{\sigma}},
\eeq
$$
\gamma_{s \sigma s}=-\frac{e_s X_{\sigma} \varphi}{(\lambda_4-\lambda_{\sigma})^2},
\quad
\gamma_{\sigma \tau \sigma}=\frac{e_{\sigma} X_{\tau} \varphi}{\lambda_{\sigma}-\lambda_{\tau}} \; (\sigma \ne \tau),
$$
where $\sigma, \tau=5, 6$, $s=2, 4$, $\gamma_{24r}$ are arbitrary, while all others $\gamma_{pqr}$ are equal to zero.

In this case the system of equations in components ${\nl{i}\xi}^j$ of skew-normal frames has the form

\medskip

$
1^{\circ}\; {\nl{1}\xi}^{\alpha}\pl_{\alpha}
{\nl{2}\xi}^{\beta}-{\nl{2}\xi}^{\alpha}\pl_{\alpha}
{\nl{1}\xi}^{\beta}=\gamma_{412}{\nl{3}\xi}^{\beta}-
\gamma_{421}{\nl{3}\xi}^{\beta}-\gamma_{411}
{\nl{4}\xi}^{\beta},
$

$
2^{\circ}\; {\nl{1}\xi}^{\alpha}\pl_{\alpha}
{\nl{3}\xi}^{\beta}-{\nl{3}\xi}^{\alpha}\pl_{\alpha}
{\nl{1}\xi}^{\beta}=\gamma_{413}{\nl{3}\xi}^{\beta}-
\gamma_{141}{\nl{1}\xi}^{\beta},
$

$
3^{\circ}\; {\nl{1}\xi}^{\alpha}\pl_{\alpha}
{\nl{4}\xi}^{\beta}-{\nl{4}\xi}^{\alpha}\pl_{\alpha}
{\nl{1}\xi}^{\beta}=-\gamma_{141}{\nl{2}\xi}^{\beta}-
\gamma_{414}{\nl{3}\xi}^{\beta}-\gamma_{241}
{\nl{1}\xi}^{\beta},
$

$
4^{\circ}\;{\nl{\sigma}\xi}^{\sigma}
\pl_{\sigma}{\nl{1}\xi}^{\beta}=
\frac{1}{2} \frac{f'_\sigma}{\lambda-f_\sigma}
{\nl{\sigma}\xi}^{\sigma}{\nl{1}\xi}^{\beta}-
\gamma_{41\sigma}{\nl{3}\xi}^{\beta},
$

$
5^{\circ}\; {\nl{2}\xi}^{\alpha}\pl_{\alpha}
{\nl{3}\xi}^{p}-{\nl{3}\xi}^{\alpha}\pl_{\alpha}
{\nl{2}\xi}^{p}=-\gamma_{142}{\nl{1}\xi}^{\beta}+
\gamma_{413}{\nl{4}\xi}^{p}+\gamma_{423}
{\nl{3}\xi}^{p},
$

$
6^{\circ}\; {\nl{3}\xi}^{\alpha}\pl_{\alpha}
{\nl{2}\xi}^{q}=-\gamma_{413}{\nl{4}\xi}^{q},
$

$
7^{\circ}\; {\nl{2}\xi}^{\alpha}\pl_{\alpha}
{\nl{4}\xi}^{\beta}-{\nl{4}\xi}^{\alpha}\pl_{\alpha}
{\nl{2}\xi}^{\beta}=-\gamma_{142}{\nl{2}\xi}^{\beta}+
\gamma_{414}{\nl{4}\xi}^{\beta}-\gamma_{242}
{\nl{1}\xi}^{\beta}+\gamma_{424}{\nl{3}\xi}^{\beta},
$

$
8^{\circ}\;{\nl{\sigma}\xi}^{\sigma}
\pl_{\sigma}{\nl{2}\xi}^{\beta}=
\frac{1}{2} \frac{f'_\sigma}{\lambda-f_\sigma}
{\nl{\sigma}\xi}^{\sigma}{\nl{2}\xi}^{\beta}-
\frac{1}{2} \frac{f'_\sigma}{(\lambda-f_\sigma)^2}
{\nl{\sigma}\xi}^{\sigma}{\nl{1}\xi}^{\beta}-
\gamma_{41\sigma}{\nl{4}\xi}^{\beta}-
\gamma_{42\sigma}{\nl{3}\xi}^{\beta},
$

$
9^{\circ}\; {\nl{3}\xi}^{\alpha}\pl_{\alpha}
{\nl{4}\xi}^{\beta}-{\nl{4}\xi}^{\alpha}\pl_{\alpha}
{\nl{3}\xi}^{\beta}=-\gamma_{143}{\nl{2}\xi}^{\beta}+
\gamma_{144}{\nl{1}\xi}^{\beta}-\gamma_{243}
{\nl{1}\xi}^{\beta},
$

$
10^{\circ}\;{\nl{\sigma}\xi}^{\sigma}
\pl_{\sigma}{\nl{3}\xi}^{\beta}=
\frac{1}{2} \frac{f'_\sigma}{\lambda-f_\sigma}
{\nl{\sigma}\xi}^{\sigma}{\nl{3}\xi}^{\beta}-
\gamma_{14\sigma}{\nl{1}\xi}^{\beta},
$

$
11^{\circ}\;{\nl{\sigma}\xi}^{\sigma}
\pl_{\sigma}{\nl{4}\xi}^{\beta}=
\frac{1}{2} \frac{f'_\sigma}{\lambda-f_\sigma}
{\nl{\sigma}\xi}^{\sigma}{\nl{4}\xi}^{\beta}-
\frac{1}{2} \frac{f'_\sigma}{(\lambda-f_\sigma)^2}
{\nl{\sigma}\xi}^{\sigma}{\nl{3}\xi}^{\beta}-
\gamma_{14\sigma}{\nl{2}\xi}^{\beta}-
\gamma_{24\sigma}{\nl{1}\xi}^{\beta},
$

$
12^{\circ}\;{\nl{\beta}\xi}^{\alpha}
\pl_{\alpha}{\nl{\sigma}\xi}^{\sigma}=0,
$

$
13^{\circ}\;{\nl{\sigma}\xi}^{\sigma}
\pl_{\sigma}{\nl{\tau}\xi}^{\tau}=
\frac{1}{2} \frac{f'_\sigma}{f_\tau-f_\sigma}
{\nl{\sigma}\xi}^{\sigma}{\nl{\tau}\xi}^{\tau},
\;\;\;(\tau\ne\sigma),
$

\noindent
where $\alpha, \beta=1, 2, 3, 4$, $p=1, 3$, $q=2, 4$, $\sigma,
\tau=5, 6$, ${\nl{\sigma}\xi}^i=P_{\sigma}(x){\delta_\sigma}^i$,
${\nl{p}\xi}^q={\nl{\alpha}\xi}^{\sigma}=0$.

\noindent
Integrating $12^{\circ}$ and $13^{\circ}$,
after coordinate transformation $x^5, x^6$ one gets
$$
{\nl{5}\xi}^5=(f_6-f_5)^{-1/2},
\quad
{\nl{6}\xi}^6=(f_5-f_6)^{-1/2}.
$$
Differentiating  $g^{\alpha \beta}={\nl{h=1}\sum^6} e_h
{\nl{h}\xi}^{\alpha} {\nl{\tilde h}\xi}^{\beta}$  with respect to $x^p$,
with using $1^{\circ}-3^{\circ}$ and
$5^{\circ}-7^{\circ}$, one gets $\pl_p g^{\alpha \beta}=0$.
Differentiating $g^{\alpha \beta}$ with respect to $x^{\sigma}$,
from equations $4^{\circ}$, $8^{\circ}$, $10^{\circ}$,
$11^{\circ}$, one obtains
$$
\pl_{\sigma} g^{pr}=-\frac{f'_\sigma}{f_\sigma-\lambda}
g^{pr}-\frac{f'_\sigma}{(f_\sigma-\lambda)^2}
{\Pi}_{\sigma}(f_{\sigma}-{\lambda})^{-1},
\quad
\pl_{\sigma} g^{pq}=-\frac{f'_\sigma}{f_\sigma-\lambda}
g^{pq},
$$
where $p, r=1, 3$. Integrating this equations, one gets
$$
g^{pr}={\Pi}_{\sigma}(f_{\sigma}-{\lambda})^{-1}({\Sigma}_1+F^{pr}),
\quad
g^{pq}=\Pi_{\sigma}(f_\sigma-\lambda)^{-1}F^{pq},
$$
where $F^{pr}, F^{pq}$ are functions of the variables $x^2, x^4$.

After same calculations for the $h$--spaces of the $[(221)1]$ type one
obtains
\bqq{24.1}
g_{ij}dx^idx^j=\Pi_{\sigma}(f_\sigma-\lambda) \lbrace 2g_{12}dx^{1}dx^{2}-e_2(\Sigma_1+\theta_1)(dx^{2})^2+
\eeq
$$
2g_{34}dx^{3}dx^{4}-e_4(\Sigma_1+\theta_2)(dx^{4})^2+G \rbrace+\sum_{\sigma} e_{\sigma} (f_\tau-f_{\sigma}) (dx^{\sigma})^2,
$$
\bqq{25}
a_{ij} dx^i dx^j=\lambda (g_{st}dx^sdx^t+G)+g_{12}(dx^{2})^2+g_{34}(dx^{4})^2+
\sum_{\sigma}f_{\sigma} g_{\sigma \sigma} (dx^{\sigma})^2+G,
\eeq
\bqq{26}
h_{ij}=a_{ij}+(\sum_{\sigma} f_\sigma+c)g_{ij},
\quad
\varphi=\frac{k f_k}{2}+\frac{1}{2}\sum_{\sigma} f_\sigma+c.
\eeq
\bqq{27}
\Sigma_1=\sum_{\sigma}(f_\sigma-\lambda)^{-1},
\eeq
\bqq{28}
G=2e_5\lbrace1+\theta_3(f-\lambda)\rbrace dx^4dx^5+(f-\lambda)g_{55} (dx^5)^2+e_6(dx^6)^2,
\eeq
where $\lambda, c$ are some constants.
Here for the $h$--spaces of the  $[(22)11]$ type  $\tau, \sigma=5, 6$ $(\tau\ne\sigma)$,
$s, t=1, 2, 3, 4$, $G=0$, $g_{st}, \theta_1, \theta_2$ are functions of the variables $x^2, x^4$,
$f_\sigma$ are functions of the variable $x^{\sigma}$.
For the $h$--spaces of the $[(221)1]$ type $\sigma=6$, $s, t=1, 2, 3, 4, 5$,
$g_{st}, \theta_1, \theta_2, \theta_3$ are functions of the variables $x^2, x^4, x^5$,
$f_6$ is a function of the variable $x^6$.

\bigskip

\noindent
{\bf Theorem 3.} {\it If a symmetric tensor $h_{ij}$ of
the characteristics $[(22)11]$, $[(221)1]$ and a function $\varphi$
satisfy in $V^6(g_{ij})$ the Eisenhart equation, then there exists a holonomic coordinate system
so that the function $\varphi$ and the tensors $g_{ij}$,
$h_{ij}$ are defined by formulas {\rm (\ref{24.1})-(\ref{28})}.}

\section{The metrics of the $h$--spaces of the $[(2211)]$, $[(22)(11)]$\\
and $[(21)(21)]$ types}

In this cases, the function $\varphi$ is constant, therefore, from (\ref{1}) the tensor  $h_{ij}$  is
covariant constant. We skip further calculations that are just integration
of the equations with respect to ${\nl{i}\xi}^j$ combined with some proper chosen coordinate
transformations. The result reads

\noindent
for the $h$--spaces of the $[(2211)]$ type
\bqq{29}
g_{ij}dx^idx^j=2g_{12}dx^1dx^2-e_2 \theta (dx^2)^2+2g_{34}dx^3dx^4+g_{rq}dx^rdx^q,
\eeq
\bqq{30}
a_{ij} dx^i dx^j=\lambda g_{ij}dx^idx^j+g_{12}(dx^2)^2,
\eeq
\bqq{31}
h_{ij}=a_{ij}+cg_{ij},
\eeq
where $r, q=5, 6$, $\lambda, c$ are some constants,
$\theta,  g_{12}, g_{34}, g_{rq}$ are functions of the variables $x^2, x^4, x^5, x^6$,

\noindent
for the $h$--spaces of the  $[(22)(11)]$ type
\bqq{32}
g_{ij}dx^idx^j=e_2\lbrace2dx^1dx^2-\theta (dx^2)^2\rbrace+
e_4\lbrace2dx^3dx^4-\omega (dx^4)^2\rbrace+
g_{\sigma \tau}dx^{\sigma}dx^{\tau},
\eeq
\bqq{33}
a_{ij} dx^i dx^j=\lambda_1 \lbrace g_{i_1j_1}dx^{i_1}dx^{j_1}+
e_2(dx^2)^2+
g_{i_2j_2}dx^{i_2}dx^{j_2}+e_4(dx^4)^2\rbrace+
\lambda_2 g_{\sigma \tau} dx^{\sigma } dx^{\tau},
\eeq
\bqq{34}
h_{ij}=a_{ij}+cg_{ij},
\eeq
where $\theta, \omega$ are functions of the variables $x^2, x^4$,
$g_{\sigma \tau}$ are functions of the variables $x^5, x^6$, $\lambda_1, \lambda_2, c$ are some constants,
here $\lambda_1 \ne \lambda_2$, $i_1, j_1=1, 2$, $i_2, j_2=3, 4$,
$\sigma, \tau=5, 6$,

\noindent
for the $h$--spaces of the $[(21)(21)]$ type
\bqq{35}
g_{ij}dx^idx^j=e_2\lbrace2dx^1dx^2-\theta (dx^2)^2\rbrace+
e_3 (dx^3)^2+
e_5 \lbrace 2dx^4dx^5-\omega (dx^5)^2\rbrace+e_6 (dx^6)^2,
\eeq
\bqq{36}
a_{ij} dx^i dx^j=\lambda_1 g_{i_1j_1}dx^{i_1}dx^{j_1}+e_2(dx^2)^2+
\lambda_2 g_{i_2j_2}dx^{i_2}dx^{j_2}+e_5(dx^5)^2,
\eeq
\bqq{37}
h_{ij}=a_{ij}+cg_{ij},
\eeq
where $\theta$ is a function of the variables $x^2, x^3$,
$\omega$ is a function of the variables $x^5, x^6$,
$\lambda_1, \lambda_2, c$ are some constants, $\lambda_1 \ne \lambda_2$,
$i_1, j_1=1, 2, 3$, $i_2, j_2=4, 5, 6$.

Now we have the following theorem.

\bigskip

\noindent
{\bf Theorem 4.} {\it If a symmetric tensor $h_{ij}$ of
the characteristics $[(2211)]$, $[(22)(11)]$, $[(21)(21)]$ and a function $\varphi$
satisfy in $V^6(g_{ij})$ equation (\ref{1}), then there exists a holonomic coordinate
system so that the function $\varphi$ and the tensors $g_{ij}$,
$h_{ij}$ are defined by formulas {\rm (\ref{29})-(\ref{37})}.}

\section{Quadratic first integrals of the geodesic\\
equation of the $h$-spaces of the $[(21\ldots1)(21\ldots1)\ldots(1\ldots1)]$ type}

For every solution $h_{ij}$ of equation (\ref{1}) there is a
quadratic first  integral of the geodesic equations (see \cite{am2},
\cite{zak})
\bqq{38}
(h_{ij}-4\varphi g_{ij})\dot{x}^i \dot{x}^j={\rm const},
\eeq
where $\dot{x}^i$ is the tangent vector to the geodesic.

Therefore, the quadratic first integrals of the geodesic equations
of the $h$-spaces of the $[(21\ldots1)(21\ldots1)\ldots(1\ldots1)]$ type are
determined by formula (\ref{38}), where the tensors $h_{ij}$, $g_{ij}$ and
the function $\varphi$ are obtained in Theorems 1-4.

\bigskip

I am grateful to professor A.V.Aminova for  constant attention to this work
and for useful discussions. The research was partially supported by the RFBR grant 01-02-17682-a
and by the INTAS grant 00-00334.

\end{document}